\documentclass[12pt,a4paper]{article}
\begin{document}

\newtheorem{Def}{Definition}
\newtheorem{Th}{Theorem}
\newtheorem{Pro}{Proposition}
\newtheorem{Le}{Lemma}
\newtheorem{Ex}{Example}
\newcommand{\be}{\begin{eqnarray}}
\newcommand{\ee}{\end{eqnarray}}
\newcommand{\bee}{\begin{eqnarray*}}
\newcommand{\eee}{\end{eqnarray*}}
\newcommand{\cl}{{\cal{C}}\!\ell}
\renewcommand{\Re}{{\cal{R}}\!\mathit{e}\,}
\renewcommand{\P}{{\cal{P}}}
\newcommand{\Q}{{\cal{Q}}}
\newcommand{\ha}{\hookrightarrow}
\newcommand{\ov}{\overline}
\newcommand{\ts}{\textstyle}
\newcommand{\im}{{\rm im}\,}
\newcommand{\nl}{\newline}
\newcommand{\vs}[1]{\ \vspace*{#1} \nl }
\renewcommand{\div}{\,{\rm div }\,}
\newcommand{\grad}{\, {\rm grad }\, }
\newcommand{\pL}{\textsl{PVec}L_{2,\cl}}
\newcommand{\BR}{{\rm\hskip 0.1pt%
                       I\hskip -2.15pt R}}

\newcommand{\BC}{{\hbox{\rm\vphantom{X}%
                       \hskip 0.25em%
                       \vrule width 0.7pt%
                       \hskip -0.35em C}}}

\noindent {\large SWANHILD BERNSTEIN}\footnote{This paper was
completed when the author was visiting the
University of Arkansas at Fayetteville, supported by a
Feodor-Lynen-fellowship of the Alexander von Humboldt foundation}
\vspace*{1ex}\\
{\huge \textbf{Monotonicity principles for singular integral
equations$\!$ in$\!$ Clifford$\!$ analysis}}
\vspace*{0.5ex}\\ \section{Introduction}
Monotonicity principles are used to get informations about
nonlinear singular integral equations. These results are based on
the famous theorem of Browder and Minty (see for example
\cite{Zei})
\begin{Th}[Browder-Minty]
Let $X$ be a real separable Banach space. If the operator $A:X\to
X^{*}$ is monotone, i.e. $\langle Au-Av\, ,\, u-v\rangle \geq 0$,
coercive, i.e. $\lim\limits_{||u||\to \infty} \frac{\langle
Au\,,\,u\rangle }{||u||}=\infty ,$ and hemicontinuous then $A$ is
surjectiv. If $A$ is moreover strictly monotone, then the solution
of $Au=b,\ b\in X^{*},$ is unique.
\end{Th}
In \cite{vW} this theorem is used by L. v. Wolfersdorf to consider
singular integral equations on the real halfline involving the
singular (complex) Hilbert operator. This theory is extended by
Askabarov in \cite{as} to the complex case.\\ We want to
investigate a special family of singular integral operators in
Clifford analysis which has in $\BR ^3$ an application to the
nonlinear magnetic field equation (\cite{frie}) considered by M.
Friedman. Electromagnetic processes are described using
quaternionic and Clifford analysis by several authors (\cite{bb},
\cite{ba}, \cite{bdsa}, \cite{g1}, \cite{g2}, \cite{gs}, \cite{h},
\cite{ks1}, \cite{ks2}, \cite{mi}, \cite{O1}, \cite{O2},
\cite{Xu1}, \cite{Xu2}). These considerations mainly based on the
operator $D+ a.$ To this subject we want to recommend the book by
V.V. Kravchenko and M.V. Shapiro (\cite{ksbook}) and the book by K.
G\"urlebeck and W. Spr\"o{\ss}ig (\cite{gs2}).\\ Singular integral
operators are investigated by S.G. Michlin and S. Pr\"o{\ss}dorf in
\cite{mp} and especially using Clifford analytical methods by A.
McIntosh, C. Li and S. Semmes in \cite{mqs} and by A. McIntosh, C.
Li and T. Qian in \cite{mlq}.\\ Properties of the Nemickii
operator, monoton operators and so one may be found in the book
\cite{Zei} by E. Zeidler.

\section{Clifford algebras} Let us denote by
$\cl_{0,h}(\BR)$ the \textit{real Clifford algebra} associated to
the Euclidean space $\BR^h$ generated by the elements
$\{e_j\}_{j=1}^h$ with
\bee e_ie_k+e_ke_i=-2\delta _{ik}.\eee
An arbitrary element of $\cl_{0,h}(\BR)$ has a representation of
the form
\bee a=\sum\limits_I a_Ie_I, \qquad a_I \in \BR, \eee
where $I$ denotes ordered l-tuples of the form $I=(i_1,i_2,\ldots
,i_l),$ with $1\leq i_1<i_2< \cdots <i_l\leq h,$ where $0\leq l\leq
h.$ Furthermore, $e_I$ stands for the product $e_{i_1}e_{i_2}\cdots
e_{i_l}.$ By convention, $e_{\emptyset}:=e_0:=1.$\\ On
$\cl_{0,h}(\BR)$ may be introduced an involution, sometimes called
main involution:
\bee  \bar{a}:= \sum\limits_I a_I\bar{e_I}, \qquad
      \bar{e_I} := (-1)^{\frac{|I|(|I|+1)}{2}}e_I.
\eee
The \textit{complexified Clifford algebra} $\cl_{0,h}(\BC)$
associated to $\BR^h$ is
\bee    \cl_{0,h}\otimes \BC, \eee
where $\BC$ denotes the complex numbers. Thus, an arbitrary element
$c$ of the complexified Clifford algebra has a representation of
form
\bee c=\sum\limits_I c_Ie_I = a+ib =\sum\limits_I a_Ie_I +i \sum\limits_I b_Ie_I,
\ c_I=a_I+ib_I \in \BC,\ a_I,b_I\in \BR. \eee
In the complexified Clifford algebra may be introduced another involution:
\bee \tilde{c} := \sum\limits_I \bar{c_I}\bar{e_I},\eee
where $\bar{c}=a_I-ib_I.$
We call $c_0e_0 = c_0 $ the scalar part of $c$ and denote by $\P $ the
operator
\bee \P :\cl _{0,h}(\BC) \to \BC:\ \P c = c_0. \eee
This operator is a projection because of $\P ^2=\P $ and we denote by $\Q$
the complementary projection $\Q=I-\P.$
Further, we denote by $\Re $ the operator
\bee \Re :\cl_{0,h}(\BC) \to \BR : \Re c= a_0. \eee
Moreover, the Clifford algebra $\cl_{0,h}(\BC)$ becomes a normed
algebra with
\bee |c|:= \left( \sum\limits_I |c_I|^2\right)^{\frac{1}{2}}= [c,c]. \eee
Here,
\bee [u,v]:= \Re \sum\limits_{I,J} \bar{u}_I\bar{e}_Iv_Je_J \eee
is a real scalar product and the Clifford algebra $\cl_{0,n}$ becomes a
real Hilbert space.

\section{Clifford analysis}
Let $G$ be a bounded domain in $\BR^n$ with Lipschitz boundary
$\partial G.$ This means that $\partial G$ may be covered by
finitely many open sets $W_k$ such that each  set $G\cap W_k$ is
represented by the inequality $x_n>g_k(x_1,x_2,\ldots , x_{n-1}), $
where $g_k$ is a Lipschitz function.\\ On such a domain, the
exterior unit normal $n(y)$ is defined for almost all $y\in
\partial G,$ and Gau{\ss}'s theorem is valid.\\ Let $u\in C^1(G)$
(taking values in a Clifford algebra) then the
\textit{Dirac operator} is defined as
\bee (Du)(x) = \sum\limits_{I}\sum\limits_{j=1}^n\frac{\partial u_I}{\partial x_j}
(x)e_je_I . \eee
Any solution $u$ of the equation $Du=0$ is called a \textit{left monogenic}
function. Analogously, if the Dirac operator is acting from the right any
solution $u$ of $uD=0$ is called a right monogenic function.\\
Moreover, we consider the \textit{disturbed Dirac operator}
\bee D_{ia}u = (D + ^{ia}\!M) = \sum\limits_{I}\sum\limits_{j=1}^n
\frac{\partial u_I}{\partial x_j}(x)e_je_I + \sum\limits_{I}\sum\limits_{j=1}^n
ia_ju_I e_je_I + ia_0 \sum\limits_{I}u_Ie_I ,\eee where
$a=a_0e_0+\sum\limits_{j=1}^na_je_j$ and $a_k\in \BR,\ k=1,2,
\ldots ,n$, thus $a$ is a paravector in the Clifford algebra $\cl _{0,n}.$\\
A fundamental solution of $D_{ia}$ is
\bee e_{ia}= e^{-i<a,x>}\{(D-^{ia_0}\!M)K_{a_0}\} ,\eee
where $<a,x>=\sum\limits_{j=1}^na_jx_j$ and
\bee K_{a_0}(x)=K_{a_0}(|x|) =\frac{1}{(2\pi )^{\frac{n}{2}}}
\left( \frac{a_0}{|x|} \right) ^{\frac{n}{2}-1}K_{\frac{n}{2}-1}(a_0|x|),
\eee
is a fundamental solution of $-\Delta +a_0^2$ (cf. \cite{ort}) and
$K_{\frac{n}{2}-1}$ denotes the modified Bessel functions, the so-called MacDonalds
functions, of order ${\frac{n}{2}-1}.$\\
With the aid of this fundamental solution we define
\bee T_{ia}u:= \int\limits_G e_{ia}(x-y)u(y) dy . \eee
Then, we get
\bee D_{ia}T_{ia}u = \left\{ \begin{array}{cl} 0 & \mbox{\ in\ } \BR^n
\backslash \bar{G} \\ u & \mbox{\ in\ } G .\end{array} \right.
\eee
As mentioned above Gau{\ss}'s formula holds and from this it is derived
\bee \int\limits_G D\, u \, dG = \int\limits_{\Gamma} n\, v\, d\Gamma .
\eee
Let $L_{2,\cl }(G)$ denote the normed space of measurable functions
$u$ from $G$ to $\cl _{0,h}$ for which the norm
\bee ||u|| = \left( \int\limits_G |u|^2 dG \right) ^{\frac{1}{2}} <\infty ,
\eee
It is a (real!) Hilbert space created by the scalar product
\bee (u,v):= \Re \int\limits_G \tilde{u}\,v\, dG =
\int\limits_G \Re (\tilde{u}\,v)\, dG .\eee

\section{The Nemyckii operator}
We want to study two types of nonsingular integral equations.
First, we require the properties of the so-called Nemyckii operator $F$ in
a Clifford analysis context. This operator is defined as
\bee (Fu)(x)=f(x,u_1(x),u_2(x), \ldots ,u_n(x)) \eee
with $u=(u_1,u_2,\ldots ,u_n).$ We make the following
assumptions:\\ (A1) Carath\'eodory condition: Let $f:G\times
\BR^n\to \cl _{0,n}(
\BC)$ be a given function, where $G$ is a nonempty set in $\BR^N,
\ n,N\geq 1.$ Moreover,
$$
x\to f(x,u) \mbox{\ is measurable on\ } G \mbox{\ for all\ }
u\in\BR^n;
$$
$$
u \to f(x,u) \mbox{\ is continuous on\ } \BR^n \mbox{\ for almost
all\ } x\in G.
$$
(A2) Growth condition: For all $(x,u)\in G\times \BR^n,$
\bee |f(x,u)|\leq a(x)+b|u|. \eee
Here, $b$ is a fixed positive number and $a\in L_2(G)$ is a real-valued
nonnegative function.\\
\begin{Pro}
Under the two assumptions (A1) and (A2), the following are
valid:\vspace*{1ex}\\
The Nemyckii operator
\bee F: L_{2,\cl}(G) \to L_{2,\cl}(G) \eee
is continuous and bounded with
\bee ||Fu||_{L_{2,\cl}} \leq const (||a||_{L_2}+ ||u||_{L_{2,\cl}}) \eee
and
\bee  (Fu,u)=\Re \int\limits_G \widetilde{f(x,u(x))}u(x) dx \quad
\mbox{\ for all\ } u\in L_{2,\cl}(G).\eee
\end{Pro}
Moreover, \vspace*{1ex}\\
\begin{tabular}{|p{7.5cm}|p{6.5cm}|} \hline
\vs{-1mm} Monotonicity of $f$: The function $f$ is monotone with respect
to $u$ i.e.
\bee [ f(x,u)-f(x,v),u-v] \geq 0 \eee
for all $ u,v\in L_{2,\cl }(G). $ \vs{-1mm}
     & \vs{7mm}  implies $F$ is monotone.  \\ \hline
\vs{-1mm} Strictly monotonicity of $f$:The function $f$ is
strictly
monotone with respect to $u$ i.e.
\bee [ f(x,u)-f(x,v),u-v] >0 \eee
for all $ u,v\in L_{2,\cl }(G). $ \vs{-1mm}
     & \vs{7mm} implies $F$ is strictly monotone. \\ \hline
\vs{-1mm} Coercivness of $f$:
\bee [f(x,u),u]\geq d|u|^2 + g(x),\eee
where $g \in L_1(G).$ \vs{-1mm}
     & \vs{1mm}  implies $F$ is coercive and
\bee (Fu,u)\geq d||u||^2 + \int\limits_G g(x) dx \eee
for all $u\in L_{2,\cl}(G).$ \vs{-1mm}\\ \hline
\vs{-1mm} Positivity of $f$: For all
$(x,u)\in G\times \BR^n$
\bee [ f(x,u),u]  \geq 0 . \eee
     & \vs{1mm}  implies
\bee (Fu,u)\geq 0 \mbox{\ for all\ }u\in L_{2,\cl}(G).
\eee \\ \hline
\end{tabular}\\
\begin{tabular}{|p{7.5cm}|p{6.5cm}|} \hline
 \vs{-1mm}Asymptotic positivity of $f$: There exists a number $R>0$ such that
\bee [ f(x,u),u] \geq 0 \eee
holds for all $(x,u)\in G\times \BR^n$ with $|u|\geq R$ and meas $G
<\infty. $ \vs{-1mm}
     & \vs{-1mm} implies
\bee (Fu,u)\geq -c \mbox{\ for all\ }u\in L_{2,\cl}(G),
\eee
where $c$ is a positive constant. \\ \hline
\vs{-1mm} Lipschitz continuity of
$f$: There is a constant $L>0$ such that
$$ |f(x,u)-f(x,v)|\leq L|u-v| .$$
     & \vs{3mm} implies $F$ is Lipschitz continuous.   \\ \hline
\end{tabular}

\section{A family of positive operators}
\subsection{The singular integral operator $D_{ia}\P T_{ia} $}
The aim of this section is to show that the operator
$D_{ia}\P T_{ia} $ is a (general) singular integral operator of
Calderon-Zygmund-type. Here, ''general'' means a singular operator plus
weakly singular operator parts.\\
For this purpose it is usefull to recall some properties of Bessel functions.
First of all the MacDonalds functions $K_p$ fulfill the recursion formula
\bee \frac{d}{dt}[ t^{-p}K_p(t)]= -t^{-p}K_{p+1}(t), \eee
second these functions are linked with the Bessel
functions of third order, the Hankel functions of second order $H_p^{(2)}$
in the following way:
\bee K_p(t)=-\frac{1}{2}\pi ie^{-\frac{1}{2}ip\pi}H_p^{(2)}\left( te^{-\frac{1}{2}
i\pi }\right) \qquad \left( -\frac{1}{2}\pi < \arg t <\pi \right) . \eee
For our singular integral operator is important the behavior $t\to 0.$
For the Hankel functions are valid
\bee H_p^{(2)}(t) \sim +i\left( \frac{2}{p} \right) ^p \frac{\Gamma (p)}{\pi} ,
\ t\to 0,\qquad (p>0) . \eee
Thus, we derive
\bee \frac{d }{d (a_0|x|)}K_{a_0}(|x|) =
\frac{1}{(2\pi )^{\frac{n}{2}}}a_0^{n-2}\cdot a_0 \frac{d}{d\,(a_0|x|)}
\left( a_0|x|\right) ^{-(\frac{n}{2}-1)}K_{\frac{n}{2}-1}(a_0|x|) = \\
= \frac{1}{(2\pi )^{\frac{n}{2}}}a_0^{n-1}(-1)\left( \frac{1}{a_0|x|}\right) ^{\frac{n}{2}-1}
K_{\frac{n}{2}}(a_0|x|).
\eee
Hence,
\bee (D-ia_0)K_{a_0}(x) = \left( \sum\limits_{j=1}^n \frac{x_je_j}{|x|} \right)
a_0\frac{d K_{a_0}(|x|)}{d (a_0 |x|)} - ia_0K_{a_0}(|x|)=\\
=\frac{-1}{(2\pi )^{\frac{n}{2}}} \left( a_0|x| \right) ^{\frac{n}{2}}
K_{\frac{n}{2}}(a_0|x|) \sum\limits_{j=1}^n \frac{x_je_j}{|x|^n}
-\frac{ia_0}{(2\pi )^{\frac{n}{2}}}\frac{1}{|x|^{n-2}}(a_0|x|) ^{\frac{n}{2}-1}
K_{\frac{n}{2}-1}(a_0|x|).\eee
and we have
$$
e_{ia}(x) = e^{-i<a,x>}\{(D-^{ia_0}\!M)K_{a_0}(x)\}=e^{-i<a,x>}\{(D-ia_0)K_{a_0}(x)\}=
$$
$$
=\frac{-e^{-i<a,x>}}{(2\pi )^{\frac{n}{2}}}\left\{ \sum \limits_{j=1}^n
\frac{x_je_j}{|x|^n}(a_0|x|)^{\frac{n}{2}}K_{\frac{n}{2}}(a_0|x|) +
\frac{ia_0}{|x|^{n-2}}(a_0|x|)^{\frac{n}{2}-1}K_{\frac{n}{2}-1}(a_0|x|)\right\} .
$$
Using the properties of modified Bessel functions we get
\bee \frac{e^{-i<a,x>}}{(2\pi )^{\frac{n}{2}}}(a_0|x|)^{\frac{n}{2}}
K_{\frac{n}{2}}(a_0|x|)=\frac{e^{-i<a,x>}}{(2\pi )^{\frac{n}{2}}}
(a_0|x|)^{\frac{n}{2}}\left(\ts{-\frac{1}{2}}\right)\pi i (-i)^{\frac{n}{2}}
H_{\frac{n}{2}}^{(2)}(-ia_0|x|) =\\
=\frac{e^{-i<a,x>}}{(2\pi )^{\frac{n}{2}}}(a_0|x|)^{\frac{n}{2}}\left(\ts{-\frac{1}{2}}
\right)\pi i (-i)^{\frac{n}{2}}i\left(\ts{\frac{2}{-ia_0|x|}}\right) ^{\frac{n}{2}}
\frac{\Gamma \left( \frac{n}{2} \right)}{\pi } + {\cal O}(|x|^{\tau }) =\\
=\frac{1}{2\pi ^{\frac{n}{2}}}\Gamma \left(\ts{\frac{n}{2}}\right)
+ {\cal O}(|x|^{\tau }) =
\frac{1}{\sigma _n} + {\cal O}(|x|^{\tau}),\ \tau >0, \mbox{\ as\ }
|x| \to 0. \eee
We simply write $e(x)$ instead of $e_0(x)$, i.e.
$$
e(x) = \frac{-1}{\sigma _n} \frac{x}{|x|^n}
$$ and thus
\bee e_{ia}-e(x) = {\cal O}(|x|^{-n+1+\tau}),\ \tau > 0,\eee
where $e(x)=e_0(x)$ and also for every component, $j=1,2,\ldots n,$ we have
\bee (e_{ia}(x))_j-(e(x))_j = {\cal O}(|x|^{-n+1+\tau}),\ \tau > 0.\eee
Now, let
\bee T_ju:= \int\limits_G (e_{ia}(x-y))_ju(y) dy ,\ j=1,2,\ldots ,n, \eee
then
\bee T_ju=\int\limits_G (e(x-y))_ju(y) dy +
\int\limits_G \{(e_{ia}(x-y))_j-(e(x-y))_j\}u(y) dy \eee
and thus (cf. \cite{mp})
$$
 \frac{\partial T_j}{\partial x_k}u =
\int\limits_G \frac{\partial }{\partial x_k}(e(x-y))_ju(y) dy +
$$
$$
\int\limits_G \frac{\partial }{\partial x_k}\left\{ (e_{ia}(x-y))_j-(e(x-y))_j
\right\} u(y) dy
-(n-2)\int\limits_{S^{n}}(x_j-y_j)\cos (r,x_k) dS^{n}
$$
$$
=\int\limits_G \frac{\partial }{\partial x_k}(e_{ia}(x-y))_ju(y) dy
-\delta _{jk}\frac{(n-2)}{n}\sigma _n\cdot u(x)
$$
and the main part of the integral operator
\bee \int\limits_G \frac{\partial }{\partial x_k}(e(x-y))_ju(y) dy
 =-\frac{1}{\sigma _n}\int\limits_G \frac{\partial }{\partial x_k}\left\{
\frac{(x_j-y_j)}{|x-y|^n}\right\} u(y) dy
\eee
is a singular integral operator of Calderon-Zygmund-type.
To sumrize, the operator $D_{ia}\P T_{ia} $ is a (general) singular
integral operator.
\subsection{Properties of the singular integral operator $D_{ia}\P T_{ia}$}
\begin{Pro}[cf. \cite{mp}]
If the symbol of a singular integral operator does not depend on
the pole then that operator is bounded in $L_2(\BR^m)$.
\end{Pro}
Therefore, the operator $D_{ia}\P T_{ia} $ is bounded in
$L_2(\BR^m)$.
\begin{Th}
The operator $D_{ia}\P T_{ia}$ is monotone in $L_2(\BR^m)$ and has a
norm less than~$1.$
\end{Th}
Proof: Because ${\cal D}_{\cl}(\BR^n)$ is dense in $L_{\cl ,2}(G)$
we get with $v\in {\cal D}_{\cl}(\BR^n)$
\bee (D_{ia}\P T_{ia}v\, ,\, v)_G = (D_{ia}\P T_{ia}v\, ,\,
D_{ia}T_{ia}v)_G\\[1ex]
=||D_{ia}\P T_{ia}v||^2_{2,G} + (D_{ia}\P T_{ia}v\, ,\, D_{ia}\Q T_{ia}v)_G
\eee
and
\bee
(D_{ia}\P T_{ia}v\, ,\, D_{ia}\Q T_{ia}v)_G=\\
= \Re \int\limits_G \widetilde{D_{ia}\P T_{ia}v} \cdot  D_{ia}\Q T_{ia}v dx =
\Re \int\limits_G -{D_{i\bar{a}}\ov{\P T_{ia}v}} \cdot  D_{ia}\Q T_{ia}v dx .\eee
On the other hand, we have
\bee \Re \int\limits_G -D\{ \ov{\P T_{ia}v} \cdot D_{ia}\Q T_{ia} v \} dx =\\
= \Re \int\limits_G -D\ov{\P T_{ia}v} \cdot D_{ia}\Q T_{ia} v dx +
\Re \int\limits_G \ov{\P T_{ia}v} \cdot (-D)D_{ia}\Q T_{ia} v dx = \\
= \Re \int\limits_G -D\ov{\P T_{ia}v} \cdot D_{ia}\Q T_{ia} v dx
+\Re \int\limits_G \ov{\P T_{ia}v} \cdot \Delta\Q T_{ia} v dx + \\
-\Re \int\limits_G \ov{\P T_{ia}v} \cdot D^{ia}\!M\Q T_{ia} v dx =\\
= \Re \int\limits_G -D\ov{\P T_{ia}v} \cdot D_{ia}\Q T_{ia} v dx
-\Re \int\limits_G \ov{\P T_{ia}v} \cdot D^{ia}\!M\Q T_{ia} v dx
\eee
Hence, we obtain
\be (D_{ia}\P T_{ia}v\, ,\, D_{ia}\Q T_{ia}v)_G= \nonumber\\
\Re \int\limits_G -D\{ \ov{\P T_{ia}v} \cdot D_{ia}\Q T_{ia} v \} dx +
\Re \int\limits_G \ov{\P T_{ia}v} \cdot D^{ia}\!M\Q T_{ia} v dx+ \nonumber\\
+\Re \int\limits_G \, ^{-i\bar{a}}\!M\ov{\P T_{ia}v} \cdot (D+^{ia}\!M)
\Q T_{ia} v \} dx = \nonumber\\
= \Re \int\limits_G -D\{ \ov{\P T_{ia}v} \cdot D_{ia}\Q T_{ia} v \} dx +
\Re \int\limits_G \ov{\P T_{ia}v} \cdot D^{ia}\!M\Q T_{ia} v \} dx+\nonumber \\
+\Re \int\limits_G \, ^{-i\bar{a}}\!M\ov{\P T_{ia}v} \cdot D\Q T_{ia} v \} dx +
\Re \int\limits_G \ov{\P T_{ia}v} \cdot (^{-i\bar{a}}\!M ^{ia}\!M) \Q T_{ia}v dx .
\label{q1}\ee
Due to $-i\bar{a}ia = \bar{a}a =|a|^2 \in \BR$ we have
 $^{-i\bar{a}}\!M^{ia}\!M=^{|a|^2}\!\!\!M$ and thus
\bee (^{-i\bar{a}}\!M ^{ia}\!M) \Q T_{ia}v = ^{|a|^2}\!\!M\Q T_{ia}v =
\Q ^{|a|^2}\!\!M T_{ia}v \in \im \Q \eee
and the last integral in (\ref{q1}) is equal to zero.
Moreover, we have
\bee ^{-\bar{a}}\!MD+ D^{a}\!M = ^{\Q{a}}\!\!MD +D^{\Q{a}}\!M
= 2\P (^{a}\!MD) \eee
and thus
\bee \Re \int\limits_G \ov{\P T_{ia}v} \cdot D^{ia}\!M\Q T_{ia} v \} dx
+\Re \int\limits_G \, ^{-i\bar{a}}\!M\ov{\P T_{ia}v} \cdot D\Q T_{ia} v \} dx =\\
=\Re \int\limits_G \ov{\P T_{ia}v} \cdot ^{-\bar{a}}\!\!MD+ D^{a}\!M \Q T_{ia} v \} dx
=\Re \int\limits\limits_G \ov{\P T_{ia}v} \cdot 2 \Q \P ^{a}\!MD T_{ia} v \} dx
= 0. \eee
To sum up, we obtain
\bee (D_{ia}\P T_{ia}v\, ,\, v)_G = \Re \int\limits_G -D\left\{ \ov{\P T_{ia}v}
 \cdot D_{ia}\Q T_{ia} v \right\} dx .
\eee
Using Gau{\ss}'s formula we get
\bee (D_{ia}\P T_{ia}v\, ,\, v)_G = \Re \int\limits_{\Gamma} -n\cdot\left\{
\ov{\P T_{ia}v} \cdot D_{ia}\Q T_{ia} v \right\} dx .
\eee
In an analogously way we derive
\bee (D_{ia}\P T_{ia}v\, ,\, v)_{\BR^n\backslash \bar{G}}=
 \Re \int\limits_{\Gamma} n\cdot\left\{
\ov{\P T_{ia}v} \cdot D_{ia}\Q T_{ia} v \right\} dx .
\eee
Due to Borel-Pompeiu formula we have
\bee 0 = ||D_{ia}T_{ia}v||^2_{\BR^n\backslash \bar{G}} =\\[1ex]
||D_{ia}\P T_{ia}v||^2_{\BR^n\backslash \bar{G}} + 2(D_{ia}\P
T_{ia}v ,D_{ia}\Q T_{ia}v)_{\BR^n\backslash \bar{G}} + ||D_{ia}\Q
T_{ia}v||^2_{\BR^n\backslash \bar{G}}.
\eee
and thus
\bee (D_{ia}\P T_{ia}v ,D_{ia}\Q T_{ia}v)_G =
-(D_{ia}\P T_{ia}v ,D_{ia}\Q T_{ia}v)_{\BR^n\backslash \bar{G}}\\[1ex]
=\ts{\frac{1}{2}}\left\{
||D_{ia}\P T_{ia}v||^2_{\BR^n\backslash \bar{G}} + ||D_{ia}\Q
T_{ia}v||^2_{\BR^n\backslash \bar{G}} \right\}
\eee
This leads to the estimation
\bee (D_{ia}\P T_{ia}v,v)_G = ||D_{ia}\P T_{ia}v||^2_G +
\frac{1}{2}\left\{
||D_{ia}\P T_{ia}v||^2_{\BR^n\backslash \bar{G}} + ||D_{ia}\Q
T_{ia}v||^2_{\BR^n\backslash \bar{G}} \right\}\geq 0
\eee
which proves the positivity of the operator $D_{ia}\P T_{ia}.$
Moreover, we get
\bee ||v||^2_G=(v,v)_G= \\[1ex]
=||D_{ia}\P T_{ia}v||^2_G +
2(D_{ia}\P T_{ia}v ,D_{ia}\Q T_{ia}v)_G + ||D_{ia}\Q T_{ia}v||^2_G
\geq ||D_{ia}\P T_{ia}v||^2_G
\eee
i.e. the norm of the operator $D_{ia}\P T_{ia}$ is not greater than
$1.$
\subsection{Nonlinear integral equations}
Let us denote by $\pL $ the subspace of all paravectors, i.e.
elements of the form $ u=\sum\limits_{j=0}^nu_je_j, $ of
$L_{2,\cl}.$ It is easy to see that $D_{ia}\P T_{ia}$ maps $\pL$
into $\pL $ and is a monotone operator. If additional $a_0=0$ then
$D_{ia}\P T_{ia}$ maps $\textsl{Vec}L_{2,\cl},$ i.e. the subspace
of all vectors $u=\sum\limits_{j=1}^nu_je_j $ into itself.
Therefore, we can use monotonicity principles to get information
about the following problem:
\bee Au=g,\qquad g\in \pL ,\\
\mbox{where}\quad Au=F(x,u)+Bu,\quad Bu=D_{ia}\P T_{ia}u.
\eee
\section{Some remarks}
\subsection{Representation formulae}
The operator $D_{ia}\P T_{ia}$ may by represented as
\bee D_{ia}\P T_{ia}u =  D_{ia} \P \left( \int\limits _G \left( e^{-i<a,x-y>}D_{-ia_0}K_{ia_0}(x-y)
 \right) u(y)dy \right) \eee
In the case $a=0$ we get simply
\bee D\P Tu = \frac{1}{2^{n}\pi ^{n/2}} D \P D\left( \int\limits _G \frac{1}{\sigma _n}\left( \frac{1}{|x-y|^{n-2}}
\right) u(y)dy \right) \eee

\subsection{An application}
In case of $n=3$ the operator $D\P T $ plays an important role in
magnetic field calculations (\cite{frie}). The general problem is
to solve the nonlinear singular integral
 equation
 \bee F(M(x),x)- \frac{1}{4\pi}\grad \div \int\limits _G \frac{M(y)}{|x-y|}dy = H_a(x) \eee
for the magnatization $M(x)=(M_1(x),M_2(x),M_3(x))$ given an
``applied field`` $H_a(x)$. $G$ is a bounded region in $\BR^3,$
which is imagined to be filled with a ferromagnetic material.
$H(x)=F(M(x),x)$ is the net field in $G$ considered as a function
of $M$; $F$ is, as a rule, an experimental function representing
the magnetic permeability of the of the ferromagnetic material,
which varies with the magnetization. this relationship is ussally
given by a single valued magnetization curve, called the $M-H$
characteristic of the magnetic material; obtained by neglecting the
hysteresis effects. The integral represents the demagnetization
field due to spatial distribution of magnetization.
\vspace*{1ex}\\
Let us denote by $H_i(x)$ the field $H(x)-H_a(x)$ induced by
$M(x).$ Then the monotonicity (=positivity) of the operator $D\P T$
leads to
\be (H_i,M)\leq 0 \label{q2} \ee
and norm not greater than one gives
\be -(H_i,M)\leq ||M||^2 \mbox{\ or\ }
\left| \left( H_i,\frac{M}{||M||} \right) \right|\leq ||M||.
\label{q3} \ee
 The relation \ref{q2} indicates that the mean angle between $H_i$ and $M$ in $G$ is not less than
$\frac{\pi}{2}$ and \ref{q3} that the mean value of the projection of $H_i$ onto
$M$ is not more than the mean value of $M$ in $G.$ \\ This
coincides with the well-known maxims of electrical engineers that
``the induced field is directed opposite of the net field, or
magnatization``, and ``the induced field is less than the
magnatization.``

\ \vspace{3cm}\\
\noindent Swanhild Bernstein\\
Institute of Applied Mathematics I\\
Department of Mathematics and Computer Science\\
University of Mining and Technology\\
D-09596 Freiberg\\
Germany\vspace*{2ex}\\
e-mail bernstein@mathe.tu-freiberg.de
\end{document}